\documentclass{amsart}

\usepackage{amsmath,arydshln,multirow}

\usepackage{hyperref}

\usepackage{cases}
\usepackage{amsmath}
\usepackage{amsfonts}
\usepackage{bm}
\usepackage{arydshln}
\usepackage{amsfonts,amsmath,amssymb,amscd,bbm,amsthm,mathrsfs,dsfont}
\usepackage{mathrsfs}
\usepackage{pb-diagram}
\usepackage{amssymb}

\usepackage[all]{xy}
\newtheorem{Thm}{Theorem}[section]

\newtheorem{Prop}[Thm]{Proposition}
\newtheorem{Lem}[Thm]{Lemma}
\newtheorem{Cor}[Thm]{Corollary}

\newtheorem{Def}[Thm]{Definition}

\newtheorem{proposition-definition}[Thm]{Proposition-Definition}

\newtheorem{conj}[Thm]{Conjecture}

\theoremstyle{remark}
\newtheorem{Rem}[Thm]{Remark}

\newtheorem{Question}[Thm]{Question}
\newtheorem{Convention}[Thm]{Conventions}
\numberwithin{equation}{section}

\begin{document}

\title{Remarks on the Boston Unramified Fontaine-Mazur Conjecture, II}


	\author{Yufan Luo}
	
	\subjclass[2020]{}
	
	\dedicatory{}
	\address{Shanghai Institute for Mathematics and Interdisciplinary Sciences (SIMIS), Shanghai 200433, China
	}
	
	\address{Research Institute of Intelligent Complex Systems, Fudan University, Shanghai
		200433, China}
	\email{yufanluo@hotmail.com}

	\maketitle
	\begin{abstract}
		In this paper, we investigate Boston’s generalization of the unramified Fontaine–Mazur conjecture for Galois representations. From a group-theoretic perspective, we first show that the conjecture can be reduced to the case of certain distinguished classes of $p$-adic analytic groups and $\mathbb{F}_{p}[[T]]$-adic analytic groups. Specifically, these are open subgroups of the groups of integral points of absolutely simple algebraic groups defined over non-Archimedean local fields. Furthermore, we provide a group-theoretic interpretation of the conjecture in terms of the virtually Golod–Shafarevich property. Finally, we establish a local–global principle and a prime-to-adjoint principle for the conjecture.
	\end{abstract}	
	
	\setcounter{tocdepth}{1}
	\tableofcontents
	
	\section{Introduction}
	\subsection{The Golod–Shafarevich Inequality}
	Let $p$ be a prime number. Let $K$ be a number field and $S$ a finite set of primes of $K$. Let $K_{S}(p)$ denote the maximal pro-$p$ extension of $K$ which is unramified outside $S$ and let $G_{K,S}(p):=\operatorname{Gal}(K_{S}(p)/K)$.
	Determining whether the Galois group $G_{K,S}(p)$ is infinite when $S$ contains no primes lying above $p$ was a long-standing open problem. This question was resolved by Golod and Shafarevich \cite{MR0161852}. For a pro-$p$ group $G$, let $d(G)$ be the generator rank of $G$ and let $r(G)$ be the relation rank of $G$. Golod and Shafarevich, with a refinement due to Vinberg \cite{MR172892}, proved that if $d(G)$ and $r(G)$ are both finite and satisfy 
	\begin{equation}\label{GSinequ}
		r(G)\leq d(G)^{2}/4, 
	\end{equation}
	then $G$ is infinite. For example, if $K$ is the rational number field $\mathbb{Q}$, then one has $d(G_{\mathbb{Q},S}(p))=r(G_{\mathbb{Q},S}(p))$. It follows that if $|S| \ge 4$, then the group $G_{\mathbb{Q},S}(p)$ is infinite.
	
	\subsection{The Boston Unramified Fontaine–Mazur Conjecture and Its Reduction} 
	The Fontaine-Mazur conjecture \cite{MR1363495}, proposed by Fontaine and Mazur, concerns $p$-adic continuous representations of the absolute Galois group of a number field. Assuming some standard conjectures in algebraic geometry, a special case of this conjecture, namely the unramified Fontaine–Mazur conjecture, is particularly appealing due to its relatively elementary formulation. It predicts that unramified 
	$p$-adic continuous representations of the absolute Galois group of a number field have finite image, cf. \cite[Conjecture 5a]{MR1363495}.
	On the other hand, as shown by Golod and Shafarevich, we have seen that unramified extensions of number fields can be infinite. The unramified Fontaine–Mazur conjecture, however, asserts that if such an extension arises as a $p$-adic Lie extension, it must be finite. 
	
	Roughly speaking, the philosophy of the unramified Fontaine-Mazur conjecture is that the eigenvalues of Frobenius elements should be roots of unity. By Chebotarev's density theorem, there is a dense subset $ D $ of the image of the $p$-adic representation such that the eigenvalues of all elements in $ D $ are roots of unity. Consequently, the image of the representation contains an open solvable subgroup, and hence is finite by class field theory. See \cite[Proposition 3.2]{MR1981910}.
	
	By a \textit{pro-$p$ ring} we mean a commutative Noetherian complete local ring with finite residue field of characteristic $ p $. Inspired by the deformation theory of Galois representations, Boston formulated a generalization of the unramified Fontaine-Mazur conjecture in the following form:
	
	\begin{conj}\cite[Conjecture 2]{MR1681626}\label{BUFM}
		Let $ K $ be a number field and $ S $ a finite set of primes of $ K $ not containing any prime above $ p $. Then every continuous homomorphism $ \rho:G_{K,S}\to {\rm GL}_{n}(A) $, where $ A $ is a pro-$p$ ring, has finite image where $ G_{K,S} $ denotes the Galois group of the maximal extension of $ K $ which is unramified outside $ S $.
	\end{conj}
	
	In the previous paper \cite{MR4976749}, we have proved that Conjecture \ref{BUFM} can be verified by restricting to the cases of $p$-adic Galois representations and $\mathbb{F}_{p}[[T]]$-adic representations. Building upon this, and in the spirit of the work of Boston \cite{MR1681626, MR2148459}, we will further reduce the conjecture to certain distinguished classes of $p$-adic analytic groups and $\mathbb{F}_{p}[[T]]$-adic analytic groups. To this end, we first introduce some terminology.
	\begin{Def}
		\begin{enumerate}
			\item Let $G $ be a profinite group and let $A$ be a pro-$p$ ring. We say that $ G$ is $A$-linear if $G$ is a closed subgroup of ${\rm GL}_{n}(A)$ for some positive integer $n$. Furthermore, $ G $ is called $ p $-adic analytic if it is $\mathbb{Z}_{p}$-linear where $\mathbb{Z}_{p}$ denotes the ring of $p$-adic integers.
			
			\item An infinite profinite group is called just-infinite if all non-trivial closed normal subgroup have finite index. A profinite group is called hereditarily just-infinite if every open subgroup is just-infinite.
			
			\item Let $G$ be a pro-$p$ group. We say that $G$ is powerful if $p$ is odd and $G/\overline{G^{p}}$ is abelian, or if $p=2$ and $G/\overline{G^{4}}$ is abelian where $\overline{G^{p}}$ denotes the closure of subgroup $G^{p}=\left\langle g^{p}~|~g\in G \right\rangle $ in $G$. Moreover, we say that $G$ is uniform if $G$ is powerful and torsion-free.
			
			\item We say that a topological group is of simple type in characteristic zero (resp. characteristic $p$) if it is isomorphic to a compact open subgroup of $\mathcal{G}(F)$, where $\mathcal{G}$ is a connected, simply connected, absolutely simple algebraic group defined over a non-Archimedean local field $F$ of characteristic zero (resp. characteristic $p$).
		\end{enumerate}
	\end{Def}

	We can now state a reduced version of the conjecture as follows.

	\begin{conj}\label{reduction}
		Let $K$ be a number field and let $p$ be a prime number.
		\begin{enumerate}
			\item If $G$ is a uniform $p$-adic analytic pro-$p$ group of simple type in characteristic zero which is also hereditarily just-infinite, then there is no pro-$p$ extension $L$ of $K$ which is ramified at finitely many primes none above $p$ such that the Galois group $\operatorname{Gal}(L/K)$ is isomorphic to $G$.
			\item If $G$ is an $\mathbb{F}_{p}[[T]]$-linear pro-$p$ group of simple type in characteristic $p$ which is also hereditarily just-infinite, then there is no unramified pro-$p$ extension $L$ of $K$ such that $\operatorname{Gal}(L/K)$ is isomorphic to $G$.
		\end{enumerate}
	\end{conj}

	For $p>2$ and $n\geq 2$, the first congruence subgroup ${\rm SL}_{n}^{1}(\mathbb{Z}_{p}):=\ker({\rm SL}_{n}(\mathbb{Z}_{p})\to {\rm SL}_{n}(\mathbb{F}_{p}))$ of the special linear group ${\rm SL}_{n}(\mathbb{Z}_{p})$ over $\mathbb{Z}_{p}$ gives an example of a group $G$ satisfying the condition in (1). Similarly, the first congruence subgroup ${\rm SL}_{n}^{1}(\mathbb{F}_{p}[[T]])$ of the special linear group ${\rm SL}_{n}(\mathbb{F}_{p}[[T]])$ over the ring of formal power series over the finite field $\mathbb{F}_p$ of order $p$, provides an example for (2). We prove the following theorem, which is the first main result of this paper and provides a further reduction of Conjecture \ref{BUFM}.

	\begin{Thm}\label{furterreduction}
		Conjecture \ref{BUFM} is equivalent to Conjecture \ref{reduction}.
	\end{Thm}
	\begin{Rem}
		\begin{enumerate}
			\item It is worth noting that, in Conjecture \ref{reduction}(2), the case of $\mathbb{F}_{p}[[T]]$-adic representations is more tractable than the $p$-adic case, as it suffices to verify the conjecture for unramified extensions; see Lemma \ref{local}.
			
			\item Note that if $F$ is a non-Archimedean local field of characteristic zero, then for any positive integer $d$ there are only finitely many isomorphism classes of simply connected simple algebraic groups defined over $F$ of dimension $d$; see \cite[Section 3.3]{MR316588}. Consequently, in the $p$-adic analytic case, Conjecture \ref{reduction}(1) admits, for each fixed dimension, only finitely many counterexamples up to commensurability. Recall that two groups are commensurable if they contain isomorphic subgroups of finite index.
			
			\item As connected, simply connected, absolutely simple algebraic groups defined over non-Archimedean local fields have already been classified, Conjecture \ref{reduction} can be verified by an inspection of the classification; cf. \cite{MR224710}.
		\end{enumerate}
	\end{Rem}

	The significance of this theorem lies in reducing the problem to a minimal counterexample to the conjecture. This is analogous to the reduction of certain problems to finite simple groups in finite group theory. The theorem precisely identifies two critical cases that suffice to establish the conjecture.
	
	\subsection{A Group-Theoretic Interpretation of Conjecture \ref{BUFM} via the Virtually Golod–Shafarevich Property} 
	To the best of the author’s knowledge, all known constructions of infinite class field towers, and more generally of infinite pro-$p$ extensions $L/K$ ramified at finitely many primes none above $p$, rely on the Golod–Shafarevich inequality \eqref{GSinequ} or its variants. It is therefore natural to ask whether every such extension arises from Golod–Shafarevich-type constructions. Inspired by Lubotzky’s conjecture on $3$-manifold groups \cite{MR707163}, Boston formulated the Virtually Golod–Shafarevich Conjecture in \cite{MR2332053}. A related question appears in \cite[Question 1]{MR1432356}. The conjecture can be stated as follows.
	
	\begin{conj}\label{VGS}
		Suppose that $L$ is an infinite pro-$p$ extension of a number field $K$ which is ramified at finitely many primes none above $p$. Then the Galois group $\operatorname{Gal}(L/K)$ contains an open pro-$p$ subgroup $G$ such that $d(G)$ and $r(G)$ are both finite and satisfy
		\begin{equation}\label{weakGSinequ}
			r(G)\leq d(G)^{2}/4.
		\end{equation}
	\end{conj}

	The second main result of this paper is the following theorem.
	
	\begin{Thm}\label{VGSimpliesBUFM}
		Assume Conjecture \ref{VGS}. Then the following assertions hold.
		\begin{enumerate}
			\item 	Conjecture \ref{BUFM} holds when $A$ is $p$-torsion-free. In particular, Conjecture \ref{reduction}(1) holds.
			\item Assume that $p$ is odd. Then Conjecture \ref{reduction}(2) holds when the corresponding algebraic group of $G$ is split.
		\end{enumerate}
		
	\end{Thm}
	\begin{Rem}
		In  \cite[Section 2]{MR2332053}, Boston claimed that Conjecture \ref{VGS} implies Conjecture \ref{BUFM} under the assumption that non-abelian free pro‑$p$ groups are not $A$-linear for any pro-$p$ ring $A$. Unfortunately, the argument contains a flaw. Specifically, Theorem 2.2 of that paper incorrectly cites a result of Zelmanov, which in fact requires the strict inequality $r(G)<d(G)^{2}/4$. The non‑strict inequality $r(G) \leq d(G)^{2}/4$ is insufficient, as counterexamples exist. For instance, the group $G = \mathbb{Z}_{p} \times \mathbb{Z}_{p}$ satisfies $d(G)=2$ and $r(G)=1$, yielding $r(G)=d(G)^2/4$, yet $G$ clearly contains no non‑abelian free pro‑$p$ subgroup. One aim of this paper is to correct this error.
	\end{Rem}

	One may ask whether strengthening the non-strict inequality (\ref{weakGSinequ}) in Conjecture \ref{VGS} to the strict inequality $r(G) < d(G)^2 /4$ would allow us to deduce Conjecture \ref{BUFM}. The answer is affirmative, and it suffices to consider the following weaker conjecture.
	
	\begin{conj}\label{strongVGS}
		Suppose that $L$ is an infinite pro-$p$ extension of a number field $K$ ramified at finitely many primes none above $p$. Then the Galois group $\operatorname{Gal}(L/K)$ contains an open Golod-Shafarevich pro-$p$ subgroup. (For the definition of a Golod–Shafarevich pro-$p$ group, see Definition \ref{defofGSgroup}.)
	\end{conj}
	
	Recall that a profinite group is said to be \textit{countably-based} if there is a countable base for the neighbourhoods of the identity. We prove the following theorem.
	
	\begin{Thm}\label{amazing}
		Assume Conjecture \ref{strongVGS}. Then the following assertions hold.
		\begin{enumerate}
			\item  If $L$ is an infinite pro-$p$ extension of a number field $K$ which is ramified at finitely many primes none above $p$, then the Galois group $\operatorname{Gal}(L/K)$ contains a non-abelian free pro-$p$ subgroup. In particular, $\operatorname{Gal}(L/K)$ contains a non-abelian free abstract subgroup.
			\item Conjecture \ref{BUFM} holds.
			\item If $L$ is an infinite pro-$p$ extension of a number field $K$ ramified at finitely many prime none above $p$, then there exists a finite extension $K'$ of $K$ and a Galois subextension $L'/K'$ of $L/K'$ such that the Galois group $\operatorname{Gal}(L'/K')$ is a just-infinite pro-$p$ group which contains an isomorphic copy of every countably-based pro-$p$ group.
		\end{enumerate}
	\end{Thm}
	
	Note that part (3) of the theorem reveals the following phenomenon: assuming Conjecture \ref{strongVGS}, the Galois group of the infinite tamely ramified pro-$p$ extension of a number field that is unramified outside a finite set of primes is an extremely large and highly complicated group, in the sense that it contains every countably-based pro-$p$ group. We prove the following result, which provides some evidence for Conjecture \ref{strongVGS}.
	
	\begin{Thm}\label{Acorollaryofmylastpaper}
		If $L$ is an infinite pro-$p$ extension of a number field $K$ ramified at finitely many primes none above $p$ such that the Galois group $\operatorname{Gal}(L/K)$ is t-linear, i.e., it is a closed subgroup of ${\rm GL}_{n}(A)$ for some commutative profinite ring $A$ and some positive integer $n$, then $\operatorname{Gal}(L/K)$ contains a non-abelian dense free abstract subgroup. In particular, \cite[Conjecture 3]{MR1432356} holds in the t-linear case.
	\end{Thm}

	\subsection{A Local–Global Principle for Conjecture \ref{BUFM}}
	In \cite[Section 2]{MR4308183}, F. Hajir, C. Maire and R. Ramakrishna gave examples of an infinite unramified Galois extension of a number field where all primes have Frobenius elements with finite order. However, it turns out that such a situation cannot occur when the corresponding Galois group is $A$-linear for any pro-$p$ ring $A$. We establish the following local–global principle for Conjecture \ref{BUFM}, which constitutes the third main result of this paper. 
	
	\begin{Thm}\label{local-globalprinciple}
		Let $p$ be a prime number, $K$ a number field and $G_{K}$ the absolute Galois group of $K$. Let $n$ be a positive integer and $\rho:G_{K}\to {\rm GL}_{n}(A)$ be a continuous homomorphism where $ A $ is a pro-$p$ ring. Assume the following conditions hold.
		\begin{enumerate}
			\item $\rho$ is unramified outside a finite set $S$ of primes of $K$.
			\item $\rho(I_{v})$ is finite for any prime $v\in S_{p}$ where $S_{p}$ denotes the set of primes of $K$ above $p$ and $I_{v}$ denotes the inertia subgroup of $G_{K}$ at $v$.
			\item $\rho(G_{v})$ is finite for all but a finite number of primes $v$ of $K$ where $G_{v}$ denotes the decomposition subgroup of $G_{K}$ at $v$.
		\end{enumerate}
		Then $\rho$ has finite image. Furthermore, there exists a finite set $T$ of primes of $K$, depending on $n$ and disjoint from $S\cup S_{p}$, with the following property: if $\rho(\operatorname{Frob}_{v})$ is unipotent for all $v\in T$, then $\rho$ has finite image where $\rho(\operatorname{Frob}_{v})$ denotes the Frobenius element corresponding to the prime $v$ in the group $\rho(G_{K})$.
	\end{Thm}
	\begin{Rem}
		\begin{enumerate}
			\item The condition in $(3)$ that $\rho(G_{v})$ is finite can be replaced by the weaker assumption that $ \rho(I_{v}) $ has finite index in $ \rho(G_{v}) $. See Lemma \ref{local}.
			\item For a pro-$p$ ring $A$ and a prime $v$ of $K$ not above $p$, one can easily construct a continuous homomorphism $G_{v}\to {\rm GL}_{n}(A)$ with infinite image. Consequently, Conjecture \ref{BUFM} cannot be simply reduced to a local problem.
		\end{enumerate}
	\end{Rem}
	
	This principle shows that the finiteness of the image of the continuous homomorphism in Conjecture \ref{BUFM} is determined by the eigenvalues at a \textit{finite set} of Frobenius elements.

	\subsection{A Prime-to-adjoint Principle for Conjecture \ref{BUFM}}
	Let $K$ be a number field and $S$ a finite set of primes of $K$. Let $A$ be a pro-$p$ ring with maximal ideal $\mathfrak{m}$ and let $\rho:G_{K,S}\to {\rm GL}_{n}(A)$ be a continuous homomorphism. When the reduction $\overline{\rho}$ of $\rho$ modulo $\mathfrak{m}$ is absolutely irreducible and $S$ contains all primes of $K$ above $p$, Boston introduced the prime-to-adjoint principle in the deformation theory of $\overline{\rho}$ in \cite{MR1079842}. We apply this method to our setting, namely the case where $S$ does not contain the primes above $p$, and obtain the following result.

	\begin{Thm}\label{primetoadjointprincple}
		Let $K$ be a number field and $\rho:G_{K,\emptyset}\to {\rm GL}_{n}(A)$ be a continuous homomorphism where $ A $ is a pro-$p$ ring with residue field $\mathbb{F}$. Assume the following conditions hold.
		\begin{enumerate}
			\item The image of $\text{Im}(\overline{\rho})$ of $\overline{\rho}$ contains a subgroup $H$ whose order is prime to $p$.
			\item The restriction $\overline{\rho}|_{H}$ of $\overline{\rho}$ on $H$ is multiplicity-free, i.e., each irreducible constituent of the associated $\mathbb{F}[H]$-module appears with multiplicity exactly one.
			\item The adjoint representation $\operatorname{ad}(\overline{\rho}|_{H})$ of $\rho|_{H}$ and $\ker(N:\operatorname{Cl}(L)/p\to \operatorname{Cl}(K')/p)$ have no irreducible subrepresentations in common as $\mathbb{F}[H]$-modules where $L$ is the splitting field of $\overline{\rho}$, $K'=L^{H}$ and $N:\operatorname{Cl}(L)/p\to \operatorname{Cl}(K')/p$ is the group homomorphism of ideal class groups induced by the norm map.
		\end{enumerate}
		Then $\rho$ has finite image. 
	\end{Thm}

	See Theorem \ref{ptatheorem} for a more general statement. Note that the previously assumed condition of “absolutely irreducible” has been replaced by the weaker requirement of being “multiplicity-free.” This principle is particularly effective when the reduction $\overline{\rho}$ of $\rho$ modulo $\mathfrak{m}$ has a relatively large image.

	\subsection{Organization and notation} 
	This paper is organized as follows. In Section \ref{section2}, we discuss a reduction of Conjecture \ref{BUFM} and prove Theorem \ref{furterreduction}. Section \ref{section3} is devoted to the virtually Golod–Shafarevich property, where Theorems \ref{VGSimpliesBUFM}, \ref{amazing}, and \ref{Acorollaryofmylastpaper} are established. Next, Section \ref{section4} investigates a local–global principle for Conjecture \ref{BUFM}, proving Theorem \ref{local-globalprinciple}. Finally, Section \ref{section5} establishes a prime-to-adjoint principle for Conjecture \ref{BUFM}, and we prove Theorem \ref{primetoadjointprincple}.

	Throughout this paper, $ p $ is a prime number, $\mathbb{F}_{p}$ is the finite field of order $p$, $\mathbb{Q}$ is the rational number field and $\mathbb{Q}_{p}$ is the field of $p$-adic numbers. We denote by $\mathbb{Z}_{p}$ the ring of $p$-adic integers and $\mathbb{F}_{p}[[T]]$ the ring of formal power series over $\mathbb{F}_{p}$.
	
	By a number field we mean a finite field extension of $\mathbb{Q}$. By a pro-$p$ extension $L$ of a number field $K$ we mean a Galois extension $L$ of $K$ such that the corresponding Galois group is a pro-$p$ group. If $K$ is a number field and $S$ a finite set of primes of $K$, then we denote by $G_{K}$ the absolute Galois group of $K$, $G_{K,S}$ the Galois group of the maximal extension of $K$ which is unramified outside $S$ and $G_{K,S}(p)$ the Galois group of the maximal pro-$p$ extension of $K$ which is unramified outside $S$.
	We also denote by $S_{p}$ the primes of $K$ above $p$. 	If $F$ is a non-Archimedean local field, then we  denote by $\mathcal{O}_{F}$ the ring of integers of $F$ and $\operatorname{char}(F)$ the characteristic of $F$. By a pro-$p$ ring we mean a commutative Noetherian complete local ring with finite residue field of characteristic $ p $.

	\section{A reduction of Conjecture \ref{BUFM}}\label{section2}
   \subsection{A result on the local Galois groups}
   	If $K$ is a number field and $v$ is a non-Archimedean prime of $K$, then we denote by $K_{v}$ the completion of $K$ at $v$, $G_{v}:=\operatorname{Gal}(\overline{K_{v}}/K_{v})$ the absolute Galois group of $K_{v}$ and $I_{v}$ the inertia subgroup of $G_{v}$.
	\begin{Lem}\label{local}
		Let $ K $ be a number field and $ v $ a non-Archimedean prime of $ K $ not above $p$. Let $F$ be a non-Archimedean local field with residue field of characteristic $p$. Then 
		\begin{enumerate}
			\item If $ \rho:G_{v}\to {\rm GL}_{n}(F) $ is a continuous homomorphism such that the group $ \rho(I_{v}) $ has finite index in $ \rho(G_{v}) $, then $\rho$ has finite image.
			\item If $F$ has positive characteristic and $ \rho:G_{v}\to {\rm GL}_{n}(F) $ is a continuous homomorphism, then $ \rho(I_{v}) $ is finite.
		\end{enumerate}
	\end{Lem}
	\begin{proof}
		Let $ \rho:G_{v}\to {\rm GL}_{n}(F)$ be a continuous homomorphism. We may assume that $\rho(G_{v})\subset {\rm GL}_{n}(\mathcal{O}_{F})$. Since $\mathcal{O}_{F}$ is a pro-$p$ ring, ${\rm GL}_{n}(\mathcal{O}_{F})$ contains a closed pro-$p$ subgroup of finite index. Since the prime $v$ does not lie above $p$, we may assume, after passing to a finite extension of $K_v$, that $\rho$ factors through the Galois group of the maximal tamely ramified extension $K_{v}^{tr}$ of $K_v$. By \cite[Theorem 7.5.3]{MR2392026}, the Galois group $ \text{Gal}(K_{v}^{tr}/K_{v}) $ is isomorphic to the profinite group topologically generated by two elements $ \tau,\sigma $ with the only relation
		\begin{align}\label{relation}
			\sigma \tau \sigma^{-1}=\tau^{q},
		\end{align}
		where the element $ \tau $ is a generator of the inertia group, and $ \sigma $ is a Frobenius lift. 
		
	  Suppose that $f:=[\rho(G_{v}):\rho(I_{v})] $ is finite. After taking another finite extension of $ K_{v} $, we can assume that $ f=1 $. Let $ L $ be the splitting field of $ \rho $, i.e., we have $ \text{Gal}(L/K_{v})=\text{Im}(\rho) $. Then the extension $ L/K_{v} $ is totally ramified which is also tamely ramified. By \cite[Corollary 4, Chapter IV]{MR554237}, $ \text{Im}(\rho) $ is procyclic. Using the relation (\ref{relation}), we see that $ \rho(I_{v}) $ is finite. This proves (1).
		
		For (2), suppose that $ F $ has positive characteristic. Then we may assume that $ \rho(\tau) $ is upper triangular after taking finite extension of $ F $. Let $ \lambda_{1},\cdots,\lambda_{n} $ denote the eigenvalues of $ \rho(\tau) $. Using the relation (\ref{relation}), we see that 
		\[ \{\lambda_{1},\cdots,\lambda_{n}\}=\{\lambda_{1}^{q},\cdots,\lambda_{n}^{q}\}. \]
		It follows that $ \lambda_{i}^{q^{n!}-1}=1 $ for all $ i $, and hence $ \lambda_{i} $ are roots of unity. Thus, $ \rho(\tau)^{m} $ is unipotent for some positive integer $ m $ which implies $ \rho(\tau)^{m} $ is a $ p$-element, and hence $ \rho(\tau) $ has finite order. This completes the proof of our lemma.
	\end{proof}
	
	\subsection{$\mathcal{O}_{F}$-analytic groups of semisimple type}  Let $F$ be a non-Archimedean local field. Let $\mathcal{O}_{F}$ denote the ring of integers of $F$ and let $\mathfrak{m}$ denote the maximal ideal of $\mathcal{O}_{F}$. For a positive integer $n$ and a set $S$, we will denote by $S^{(n)}$ the $n$-th Cartesian power of $S$. For the convenience of the reader, we recall some basic definitions on the theory of $\mathcal{O}_{F}$-analytic groups.
	\begin{Def}
		\begin{enumerate}
			\item \cite[Definition 13.1]{MR1720368} Let $U$ be an open subset of $\mathcal{O}_{F}^{(n)}$. A map $f:U\to \mathcal{O}_{F}^{(m)}$ is analytic if for each $x\in U$ there exists a positive integer $N$ such that $x+(\mathfrak{m}^{N})^{(n)}\subset U$ and a tuple of power series $\mathbf{F}\in \mathcal{O}_{F}[[X_{1},\cdots,X_{n}]]^{(m)}$ such that $f(y)=\mathbf{F}(x+y)$ for all $y\in (\mathfrak{m}^{N})^{(n)}$.
			\item \cite[Definition 13.5]{MR1720368} Let $X$ be a topological space and $U$ a non-empty open subset of $X$. A triple $(U,\phi,n)$ is an $\mathcal{O}_{F}$-chart on $X$ if $\phi$ is a homeomorphism from $U$ onto an open subset of $\mathcal{O}_{F}^{(n)}$ for some $n\in \mathbb{N}$. Two charts $(U,\phi,n)$ and $(V,\psi,m)$ on $X$ are compatible if the maps $\psi\circ \phi^{-1}|_{\phi(U\cap V)}$ and $\phi \circ \psi^{-1}|_{\psi(U\cap V)}$ are $\mathcal{O}_{F}$-analytic functions on $\phi(U\cap V)$ and $\psi(U\cap V)$ respectively. A $\mathcal{O}_{F}$-atlas on $X$ is a set of pairwise compactible $\mathcal{O}_{F}$-charts that cover $X$. Two atlases $A$ and $B$ on $X$ are compatible if every chart in $A$ is compatible with every chart in $B$.
			
			An $\mathcal{O}_{F}$-analytic structure over $\mathcal{O}_{F}$ on a topological space $X$ is an equivalence class of compatible atlases. An $\mathcal{O}_{F}$-analytic manifold is a topological space $X$ endowed with an $\mathcal{O}_{F}$-analytic structure. 
			
			\item \cite[Definition 13.8]{MR1720368} A $\mathcal{O}_{F}$-analytic group is a topological group $G$ that is an $\mathcal{O}_{F}$-analytic manifold in such a way that the multiplication map $m:G\times G\to G,(x,y)\mapsto xy$ and the inversion map $\tau:G\to G,x\mapsto x^{-1}$ are analytic maps.
			
			\item \cite[Definition 13.13]{MR1720368} An $\mathcal{O}_{F}$-standard group $G$ of level $n$ and dimension $d$ is an $\mathcal{O}_{F}$-analytic group such that
			\begin{enumerate}
				\item the $\mathcal{O}_{F}$-analytic structure on $G$ can be defined by a global atlas $\{(G,\psi,d)\}$ where $\psi=(\psi_{1},\cdots,\psi_{d})$ is a homeomorphism onto $(\mathfrak{m}^{n})^{(d)}$ and $\psi(1)=0$;
				\item for $j=1,\cdots,d$, there exists formal power series $F_{j}(\mathbf{X},\mathbf{Y})\in \mathcal{O}_{F}[[\mathbf{X},\mathbf{Y}]]$, without constant term, such that $\psi(xy)=\mathbf{F}(\psi(x),\psi(y))$ for all $x,y\in G$ where $\mathbf{X}=(X_{1},\cdots,X_{d})$ and $\mathbf{Y}=(Y_{1},\cdots,Y_{d})$.
			\end{enumerate}
			We use the convention that $\mathcal{O}_{F}$-standard group means $\mathcal{O}_{F}$-standard group of level $1$.
			\item \cite[Definition 13.28]{MR1720368} Let $G$ be an $ \mathcal{O}_{F}$-standard group of positive dimension $d$. We identify $G$ with $\mathfrak{m}^{(d)}$. We say that $G$ is $\mathcal{O}_{F}$-perfect if $[G,G]=(\mathfrak{m}^{2})^{(d)}$.
		\end{enumerate}
 
	\end{Def}
	For more details, we refer the reader to \cite[Chapter 13]{MR1720368}.

	\begin{Def}
		We say that a topological group is of \textit{semisimple type} in characteristic zero (resp. characteristic $p$) if it is isomorphic to a compact open subgroup $\Gamma$ of $\mathcal{G}(F)$, where $\mathcal{G}$ is a connected, simply connected, semisimple algebraic group defined over a non-Archimedean local field $F$ of characteristic zero (resp. characteristic $p$). Furthermore, we say that a topological group of semisimple type is of simple type if the corresponding algebraic group is absolutely simple.
	\end{Def}
	\begin{Rem}\label{rem}
	The subgroup $\Gamma$ is Zariski-dense in $\mathcal{G}$. Furthermore, every topological group of semisimple type is a finitely presented profinite group and admits the structure of an $\mathcal{O}_{F}$-analytic group. For these facts we refer to \cite[Section 2]{MR2363421}.
	\end{Rem}
	
	\begin{Prop}\label{splitcaseisperfect}
	Let $\mathcal{G}$ be a connected, simply connected, absolutely simple, split algebraic group defined over a non-Archimedean local field $F$ of characteristic $p$. Then the group $\ker(\mathcal{G}(\mathcal{O}_{F})\to \mathcal{G}(\mathcal{O}_{F}/\mathfrak{m}))$ is $\mathcal{O}_{F}$-perfect where $\mathfrak{m}$ denotes the maximal ideal of $\mathcal{O}_{F}$. 
	\end{Prop}
	\begin{proof}
	It is \cite[Exercise 11, Chapter 13]{MR1720368} or \cite[Remark 5.3]{MR1264349}.
	\end{proof}
	
	We have the following important theorems of Pink, which describe the structure of compact linear groups over local fields.
	
	\begin{Thm}[Pink]\label{pink}
		Let $F$ be a non-Archimedean local field and let $\Gamma$ be a finitely generated compact subgroup of ${\rm GL}_{n}(F)$. Then there exist closed normal subgroups $\Gamma_{2}\subset \Gamma_{1}\subset \Gamma$ such that:
		\begin{enumerate}
			\item $\Gamma/\Gamma_{1}$ is finite.
			\item $\Gamma_{1}/\Gamma_{2}$ is trivial or an $\mathcal{O}_{F}$-analytic group of semisimple type in characteristic $\operatorname{char}(F)$.
			\item $\Gamma_{2}$ is a solvable group.
		\end{enumerate}
	\end{Thm}
	\begin{proof}
 It is \cite[Corollary 0.5]{MR1637068}. See also \cite[Theorem 2.2]{MR2363421}.
	\end{proof}
	
	\begin{Thm}[Pink]\label{justinfinite}
	 Let $\Gamma$ be a finitely generated compact open subgroup of $\mathcal{G}(F)$ where $\mathcal{G}$ is a connected, simply connected, absolutely simple algebraic group defined over a non-Archimedean local field $F$. Then $\Gamma$ contains an open hereditarily just-infinite $\mathcal{O}_{F}$-analytic group of semisimple type.
	\end{Thm}
	\begin{proof}
	By Remark \ref{rem}, $\Gamma$ is Zariski-dense in $\mathcal{G}(F)$. Then it follows from \cite[Theorem 2.4]{MR2363421}.
	\end{proof}
	
	\subsection{Proof of Theorem \ref{furterreduction}}
		It is clear that Conjecture \ref{BUFM} implies Conjecture \ref{reduction}. For the other direction, it suffices to prove Conjecture \ref{BUFM} for $p$-adic representations and $\mathbb{F}_{p}[[T]]$-adic representations by \cite[Theorem 1.2]{MR4976749}.

		Now, let $\rho:G_{K,S}\to {\rm GL}_{n}(\mathcal{O}_{F})$ be a continuous homomorphism where $F$ is a non-Archimedean local field with residue field of characteristic $p$. By Lemma \ref{local}(2), when $F$ has positive characteristic, we may assume that $S=\emptyset$ after taking a finite extension of $K$. Since ${\rm GL}_{n}(\mathcal{O}_{F})$ contains a closed pro-$p$ subgroup of finite index, we may further assume that the image $\text{Im}(\rho)$ of $\rho$ is a pro-$p$ group after replacing $K$ by another finite extension. Suppose that $\Gamma:=\text{Im}(\rho)$ is infinite. By Theorem \ref{pink}, there exist closed normal subgroups $\Gamma_{2}\subset \Gamma_{1}\subset \Gamma$ such that:
		\begin{enumerate}
			\item $\Gamma/\Gamma_{1}$ is finite.
			\item $\Gamma_{1}/\Gamma_{2}$ is trivial or an $\mathcal{O}_{F}$-analytic group of semisimple type in characteristic $\operatorname{char}(F)$.
			\item $\Gamma_{2}$ is a solvable group.
		\end{enumerate}
		After taking a finite extension of $K$, we can assume that $\text{Im}(\rho)=\Gamma_{1}$. Since $\Gamma_{1}$ is infinite and $\Gamma_{2}$ is solvable, $\Gamma_{1}/\Gamma_{2}$ must be an infinite $\mathcal{O}_{F}$-analytic group of semisimple type by \cite[Theorem 2.4]{MR4976749}. The corresponding semisimple algebraic group is a finite product of simple components. Using the argument in the proof of \cite[Theorem 2.3]{MR2363421}, we may assume further that $G:=\Gamma_{1}/\Gamma_{2}$ is of simple type in characteristic $\operatorname{char}(F)$. By Theorem \ref{justinfinite}, $G$ contains an open hereditarily just-infinite $\mathcal{O}_{F}$-analytic subgroup of semisimple type. Descending again to open subgroups and simple components, we may assume that $G$ an $\mathcal{O}_{F}$-analytic pro-$p$ group of simple type which is also hereditarily just-infinite.
		
		If $\operatorname{char}(F)=0$, then we may assume that $G$ is a uniform $p$-adic analytic pro-$p$ group of simple type in characteristic zero which is also hereditarily just-infinite using the fact that any $p$-adic analytic group contains an open uniform subgroup by \cite[Interlude A]{MR1720368}. This contradicts Conjecture \ref{reduction}, and hence $\rho$ has finite image.
		
		If $\operatorname{char}(F)=p$, then $G$ is $p$-adic analytic or is a closed subgroup of ${\rm GL}_{N}(\mathbb{F}_{p}[[T]])$ for some integer $N$ by \cite[Theorem 1.7]{MR1935507}. Since $G$ is infinite, it can not be $p$-adic analytic by \cite[Theorem 1.8]{MR1935507} and \cite[Theorem 2.4]{MR4976749}. We conclude that $G$ is an $\mathbb{F}_{p}[[T]]$-linear pro-$p$ group of semisimple type which is also hereditarily just-infinite. This contradicts Conjecture \ref{reduction}. Thus, $\rho$ has finite image. The proof is finished.

	\section{The virtually Golod–Shafarevich property}\label{section3}
	\subsection{Proof of Theorem \ref{VGSimpliesBUFM}}
	
	Let $G$ be a pro-$p$ group. We define the \textit{dimension subgroups} of $G$ recursively as follows.
	\begin{Def}\cite[Definition 11.1]{MR1720368}
		$D_{1}(G)=G$, and for $i>1$
		\[ D_{i}(G):=D_{\lceil i/p \rceil }(G)^{p}\cdot \prod_{j+k=i}[D_{j},D_{k}] \]
		where $\lceil i/p \rceil $ is the least integer $r$ such that $pr\geq i$.
	\end{Def}
	\begin{Lem}\cite[Theorem D1, Interlude D]{MR1720368}\label{interlude}
		Let $G$ be a finitely generated pro-$p$ group with dimensional subgroup series $\{D_{i}(G)\}$. Suppose that $d(G)>1$ and $[G:D_{i}(G)]=p^{s_{i}}$. If $\limsup_{i\to \infty}s_{i}^{1/i}\leq 1$, then
		\[ r(G)\geq d(G)^{2}/4, \]
		and the inequality is strict unless $d(G)=2$ and $r(G)=1$.
	\end{Lem}
	
	Inspired by \cite[Theorem 4.6.4]{MR1978431}, we introduce the following definition.
	
	\begin{Def}
	Let $G$ be a finitely presented pro-$p$ group. If there exists $\epsilon>0$ such that
	\[ a_{n}(G)\leq n^{(\log_{2} n)^{2-\epsilon}} \]
	for all sufficiently large $n$ where $a_{n}(G)$ denotes the number of open subgroups of index $n$ in $G$, then we say that $G$ is admissible.
	\end{Def}

	\begin{Lem}\label{admissibleiscommen}
		Let $G$ be a finitely presented pro-$p$ group and $H$ be an open subgroup of $G$. If $G$ is admissible, then so is $H$.
	\end{Lem}
    \begin{proof}
    Let $H \leq G$ be an open subgroup of finite index $m = [G:H]$. By \cite[Lemma 2]{MR1109625}, $H$ is also finitely presented. If $G$ is admissible, then for large $n$ we have 
    \[  a_{n}(H)\leq a_{nm}(G)\leq (mn)^{(\log_{2}(mn))^{2-\epsilon}}\]
    for some $\epsilon>0$. For large $n$, we can write
    \[
    \log_2(mn) = \log_2 n + O(1).
    \]
    It follows that 
    \[
    a_n(H) \le (mn)^{(\log_2 n + O(1))^{2-\epsilon}} 
    \sim n^{(\log_2 n)^{2-\epsilon}} \cdot m^{(\log_2 n)^{2-\epsilon}}.
    \]
    The factor $m^{(\log_2 n)^{2-\epsilon}}$ grows slower than any $n^\delta$ for fixed $\delta>0$. Hence, by possibly choosing a slightly smaller $\epsilon'>0$, we see that $H$ is also admissible.
    \end{proof}

	\begin{Thm}\label{nice}
		Let $G$ be a finitely presented pro-$p$ group with $d(G)>1$. Suppose that $G$ does not satisfy $d(G)=2$ and $r(G)=1$. If $G$ is admissible, then we have
		\[ r(G)>d(G)^{2}/4. \]
	\end{Thm}  
	\begin{proof}
		The same argument as in the proof of \cite[Theorem 4.6.4]{MR1978431} shows that $\limsup_{i\to \infty}s_{i}^{1/i}\leq 1$. The claim then follows from Lemma \ref{interlude}.
	\end{proof}
	
	We thus obtain the desired theorem.
	
	\begin{Thm}\label{p-adicanalyticareGS}
		Let $G$ be a finitely presented pro-$p$ group with $d(G)>1$. Suppose that $G$ does not satisfy $d(G)=2$ and $r(G)=1$. If $G$ is an open subgroup of a $p$-adic analytic group or an $\mathcal{O}_{F}$-perfect group for some non-Archimedean local field $F$ of characteristic $p$, then 
		\[ r(G)>d(G)^{2}/4. \]
	\end{Thm} 
	\begin{proof}
	It's well known that $p$-adic analytic groups and $\mathcal{O}_{F}$-perfect groups are admissible, cf. \cite[Section 4.1  and Section 4.4]{MR1978431}. Then our theorem follows from Lemma \ref{admissibleiscommen} and Theorem \ref{nice}.
	\end{proof}
	
	The following observation is trivial but useful.
	\begin{Lem}\label{simpleobst}
		Let $G$ be a finitely presented pro-$p$ group. If the abelianization $G^{ab}$ of $G$ is finite, then $r(G)\geq d(G)$.
	\end{Lem}
	
	As a consequence, we obtain
	\begin{Lem}\label{smallobs}
		Suppose that $L$ is an infinite pro-$p$ extension of a number field $K$ which is ramified at finitely many primes none above $p$. Then the Galois group $G:=\operatorname{Gal}(L/K)$ satisfies $d(G)>1$ and $r(G)\geq d(G)$. In particular, $G$ does not satisfy $d(G)=2$ and $r(G)=1$.
	\end{Lem}
	\begin{proof}
		The group $G$ cannot be abelian and the abelianization of $G$ is finite by class field theory. Thus, we have $d(G)>1$ and $r(G)\geq d(G)$ by Lemma \ref{simpleobst}.
	\end{proof}
	
	We say that two groups are commensurable if they contain isomorphic subgroups of finite index. Now we can prove Theorem \ref{VGSimpliesBUFM} as follows.
	\begin{proof}[Proof of Theorem \ref{VGSimpliesBUFM}]
		For (1), by \cite[Theorem 1.2(1)]{MR4976749}, it suffices to prove Conjecture \ref{BUFM} for $A=\mathbb{Z}_{p}$. Let $\rho:G_{K,S}\to {\rm GL}_{n}(\mathbb{Z}_{p})$ be a such continuous homomorphism. After taking a finite extension of $K$, we may assume that the image $G$ of $\rho$ is a pro-$p$ group. Suppose that $G$ is infinite.  By Lemma \ref{smallobs} and Theorem \ref{p-adicanalyticareGS}, we have $r(U)>d(U)^{2}/4$ for any open subgroup $U$ of $G$. This contradicts Conjecture \ref{VGS}. We conclude that $G$ is finite, i.e., $\rho$ has finite image. This proves (1).
		
		For (2), let $L$ be an infinite unramified pro-$p$ extension of a number field $K$ such that the corresponding Galois group $G:=\operatorname{Gal}(L/K)$ is a compact open subgroup of $\mathcal{G}(F)$ where $\mathcal{G}$ is a connected, simply connected, absolutely simple, split algebraic group defined over a non-Archimedean local field $F$ of characteristic $p$. By Proposition \ref{splitcaseisperfect}, $G$ is commensurable to an $\mathcal{O}_{F}$-perfect pro-$p $ group. After taking a finite extension of $K$, we can assume that $G$ is an open subgroup of an $\mathcal{O}_{F}$-perfect pro-$p$ group. By Lemma \ref{smallobs} and Theorem \ref{p-adicanalyticareGS}, we have $r(U)>d(U)^{2}/4$ for any open subgroup $U$ of $G$. This contradicts Conjecture \ref{VGS}. This completes the proof of our theorem.
	\end{proof}
	
	Combining this theorem with Theorem \ref{furterreduction}, we immediately obtain the following corollary.
	
	\begin{Cor}\label{nextjiay}
		Suppose that $p$ is odd. The following assertions are equivalent.
		\begin{enumerate}
			\item Conjecture \ref{VGS} implies Conjecture \ref{BUFM}.
			\item  Conjecture \ref{VGS} implies Conjecture \ref{reduction}(2) when the corresponding algebraic group of $G$ is not split.
		\end{enumerate}
	\end{Cor}
	
	\subsection{Golod-Shafarevich pro-$ p$ groups}Let $G$ be a finitely generated pro-$p$ group and let $I(G)$ be the augmentation ideal of the completed algebra $\mathbb{F}_{p}[[G]]$. For each positive integer $n$, let $G_{n}=\{g\in G~|~g-1\in I(G)^{n}\}$. The series $\{G_{n}\}$ is called the \textit{Zassenhaus filtration} of $G$. 	For more details, we refer the reader to \cite[Section 7.4]{MR1930372}. We introduce the definition of Golod-Shafarevich pro-$p$ groups, which coincides with the one in \cite[Section 3.1]{MR2949205}.

	\begin{Def}\label{defofGSgroup}
		Let $G$ be a finitely generated pro-$p$ group with $d$ generators and let 
		\begin{equation}\label{presentation}
			1\to R\to F\to G\to 1
		\end{equation}
		be a presentation of $G$ where $F$ is a finitely generated free pro-$p$ group of rank $d$. The \textit{level} of $r\in R$ is the natural number $m$ such that $r\in F_{m}\backslash F_{m+1}$ where $\{F_{m}\}$ is the Zassenhaus filtration. Choose a generating set $\{\rho_{i}\}_{i\in I}$ for $R$ as a normal subgroup of $F$, and let $r_{k}$ denote the number of these relations which have level $k$. Put
		\[ P_{R}(t):=1-dt+\sum_{k\geq 2}r_{k}t^{k}.\]
		If there exists a real number $t_{0}\in (0,1)$ such that $P_{R}(t_{0})<0$, then we say that the pro-$ p$ presentation (\ref{presentation}) of $G$ satisfies the \textit{Golod-Shafarevich condition}. A pro-$p$ group $G$ is called a \textit{Golod-Shafarevich group} if has a presentation satisfying the Golod-Shafarevich condition.
	\end{Def}
	
	For example, if $G$ is a finitely presented pro-$p$ group such that $r(G)<d(G)^{2}/4$ and $d(G)>1$, then $G$ is Golod-Shafarevich, cf. \cite[Theorem 3.4]{MR2949205}. The Golod-Shafarevich pro-$p$ groups satisfy the following properties.
	
	\begin{Prop}\label{propertiesofGSgroup}
		\begin{enumerate}
			\item Every Golod-Shafarevich pro-$p$ group is infinite.
			\item Every Golod-Shafarevich pro-$p$ group contains a non-abelian free pro-$p$ subgroup.
			\item If $G$ is a Golod-Shafarevich pro-$p$ group, then it has a quotient $G'$ such that $G'$ is Golod-Shafarevich which contains a dense subset of torsions. 
			\item Every Golod-Shafarevich pro-$p$ group is not just-infinite.
		\end{enumerate}
	\end{Prop}
	\begin{proof}
		The statement (1) is due to Vinberg \cite{MR172892}, see also \cite[Theorem 7.20]{MR1930372} and \cite[Theorem 3.1]{MR2949205}. The statement (2) is due to Zelmanov in  \cite{MR1765122}. The statement (3) is due to Wilson, see \cite[Theorem A']{MR1109625} and \cite[Theorem 6.2]{MR2949205}. The statement (4) follows from \cite[Observation 6.3]{MR2949205}.
	\end{proof}

  \subsection{Proof of Theorem \ref{amazing}}
	We say that a profinite group $ G$ has a property $\mathcal{P}$ virtually if $G$ has an open subgroup with property $\mathcal{P}$. A branch profinite group is a certain kind of subgroup of the automorphism group of a locally finite, rooted tree. For more details, see \cite{MR1765119}. We only need the following fact.

	\begin{Thm}\label{Grigo}
		\begin{enumerate}
			\item Let $G$ be a profinite just infinite group. Then either $G$ is a branch group, or $G$ contains an open normal subgroup which is isomorphic to the direct product of a finite number of copies of some hereditarily just infinite profinite group.
			\item Branch pro-$p$ groups, which are not virtually torsion-free, contain a copy of every countably based pro-$p$ group.
		\end{enumerate}
	\end{Thm}
	\begin{proof}
		The assertion (1) is \cite[Theorem 3]{MR1765119}, and the assertion (2) is \cite[Theorem 2]{MR1754662}.
	\end{proof}
	
	Now we can prove Theorem \ref{amazing} as follows.
		\begin{enumerate}
			\item It follows from Proposition \ref{propertiesofGSgroup}(2).
			\item By \cite[Theorem 1.2]{MR4976749}, it suffices to prove Conjecture \ref{BUFM} for $A=\mathcal{O}_{F}$ where $F$ is a non-Archimedean local field with residue field of characteristic $p$. Suppose that $\rho$ has infinite image $G$. Then $G$ must contain a non-abelian free pro-$p$ subgroup by (1). But this is impossible by \cite[Theorem 1.1]{MR1684856}. Therefore, $\rho$ has finite image.
			
			\item Suppose that $L$ is an infinite pro-$p$ extension of a number field $K$ ramified at finitely many prime none above $p$. By Conjecture \ref{strongVGS}, there is a finite extension $K'$ of $K$ such that $\operatorname{Gal}(L/K')$ is Golod-Shafarevich. By Proposition \ref{propertiesofGSgroup}, $\operatorname{Gal}(L/K')$ has an infinite quotient $\Gamma$ which contains a dense subset of torsions. By \cite[Proposition 3(b)]{MR1765119}, $\Gamma$ can be mapped onto a just-infinite pro-$p$ group $\Gamma'$. By Conjecture \ref{strongVGS}, Proposition \ref{propertiesofGSgroup}(4) and Theorem \ref{Grigo}(1), $\Gamma'$ must be a profinite branch group. Since $\Gamma$ contains a dense subset of torsions, $\Gamma'$ is not virtually torsion-free. By Theorem \ref{Grigo}(2), $\Gamma'$ contains an isomorphic copy of every countably-based pro-$p$ group. Let $L'/K'$ be the Galois subextension of $L/K'$ corresponding to $\Gamma'$. Then $\operatorname{Gal}(L'/K')$ is a just-infinite pro-$p$ group which contains an isomorphic copy of every countably-based pro-$p$ group. This completes the proof of our theorem.
		\end{enumerate}

	\begin{Rem}
		Theorem \ref{amazing} remains valid if the Golod-Shafarevich pro-$p$ subgroup in the statement of Conjecture \ref{strongVGS} is replaced by the generalized Golod-Shafarevich pro-$p$ group defined in \cite[Section 4.2]{MR2949205}.
	\end{Rem}

	\subsection{Proof of Theorem \ref{Acorollaryofmylastpaper}}
By \cite[Theorem 11.5]{MR1930372}, the Galois group $\operatorname{Gal}(L/K)$ is finitely generated. By \cite[Theorem 1.1]{Luo22012026}, the Galois group $\operatorname{Gal}(L/K)$ contains either a solvable subgroup of finite index or a non-abelian dense free abstract subgroup. Then our claim follows from \cite[Theorem 2.4]{MR4976749}.

	\subsection{Further questions}
	Let $G$ be a pro-$p$ group with finite generator rank $d(G)$ and finite relation rank $r(G)$. If $G$ satisfies $r(G) < d(G)^2/4$ 
	and has $d(G) \ge 2$, then $G$ contains a non-abelian free pro-$p$ subgroup by Proposition \ref{propertiesofGSgroup}(2). In the borderline case $r(G) = d(G)^2/4$, the group is still infinite, 
	but it may fail to contain such a non-abelian free pro-$p$ subgroup. For example, the group $G = \mathbb{Z}_p \times \mathbb{Z}_p$ has $d(G)=2$ and $r(G)=1$ but $G$ does not contain a non-abelian free pro-$p$ subgroup. This leads to the following question.
	
	\begin{Question}
		Let $G$ be a finitely generated pro-$p$ group satisfying $r(G)=d(G)^{2}/4$ where $d(G)>2$. Does $G$ contain a non-abelian free pro-$p$ subgroup?
	\end{Question}
	
	Furthermore, Corollary \ref{nextjiay} naturally leads to the following question:
	
	\begin{Question}
		Let $G$ be a finitely presented pro-$p$ group with $d(G)>1$. Suppose that $G$ does not satisfy $d(G)=2$ and $r(G)=1$. If $G$ is commensurable to $\mathcal{G}(\mathcal{O}_{F})$ where $\mathcal{G}$ is a connected, simply connected, absolutely simple algebraic group defined over a non-Archimedean local field $F$ of characteristic $p$, then is it true that
	\[ r(G)>d(G)^{2}/4?\]
	\end{Question}
	
   As proved earlier, the answer to this question is affirmative if $\mathcal{G}$ is split over $F $. Thus, it suffices to consider the non‑split case. If the answer remains affirmative in this setting, then the argument presented above would indeed show that Conjecture \ref{VGS} implies Conjecture \ref{BUFM}. We may further ask:
	
		\begin{Question}
		Is Conjecture \ref{VGS} equivalent to Conjecture \ref{strongVGS}?
	\end{Question}
	
	\section{A local-global principle}\label{section4}
	\subsection{Torsion elements in a $A$-linear group}
	Let $ A $ be a pro-$p$ ring. We say that $A$ is a \emph{pro-$p$ domain} if it is an integral domain. For a pro-$p$ domain $(A,\mathfrak{m}_{A})$, we put $A_{0}:=\mathbb{Z}_{p}$ when $\text{char}(A)=0$ and $A_{0}:=\mathbb{F}_{p}[[T]]$ when $\text{char}(A)=p$. The field of fractions of the integral domain $ A_{0}$ is denoted by $\text{Frac}(A_{0})$. We prove that:
	
	\begin{Thm}\label{niceobsection}
		Let $G$ be a closed subgroup of $ {\rm GL}_{n}(A)$ where $A$ is a pro-$p$ ring. If $G$ contains a dense subset $D$ of torsions, then it contains a solvable subgroup of finite index.
	\end{Thm}
	\begin{proof}
		First of all, we show that the claim holds for $A=\mathcal{O}_{F}$ where $F$ is some non-Archimedean local field. By \cite[Lemma 3.4]{MR4976749}, there exists a positive integer $ m $ such that $ x $ satisfies the polynomial $ t^{m}-1=0 $ for any $ x\in D $. Hence the minimal polynomial of $ x $ must divide this polynomial for any $ x\in D $. In other words, $G$ contains a dense subset $ D $ such that the roots of the characteristic polynomial $ P_{x}(t)=\det(1-t\cdot x) $ of $ x $ are roots of unity for any $x\in D$, and these roots of unity have bounded order uniformly as $x $ varies over $ D $. 
		Since the degree of the characteristic polynomial $ P_{x}(t) $ for some $ x\in D $ is fixed and independent of $x\in D  $, there are only finitely many possibilities for the characteristic polynomial $ P_{x}(t) $ as $ x $ varies over $ D $. Since $ D $ is dense in $ G $, we see that there are also only finitely many possibilities for the characteristic polynomial of any element in $ G$, and the characteristic polynomial map $ G\to F[t] $ sending any element of $ G $ to the characteristic polynomial of the element is locally constant. It follows that the set $U= \{g\in G~|~\det(1-tg)=(1-t)^{n}\} $ is open in $ G $. Therefore, $G$ contains a solvable subgroup of finite index. 
		
		In general, let $A$ be a pro-$p$ ring and let $A^{\text{red}}$ denote the quotient of $A$ by the nilradical of $A$. By \cite[Lemma 3.8]{MR4976749}, the kernel of the group homomorphism $ {\rm GL}_{n}(A)\twoheadrightarrow {\rm GL}_{n}(A^{\text{red}}) $ induced by the ring homomorphism $ A\twoheadrightarrow A^{\text{red}} $ is a nilpotent group. Thus, after replacing $A$ by $A^{\text{red}}$, we may assume that $A$ is reduced. Then we consider the injective group homomorphism
		\[{\rm GL}_{n}(A)\hookrightarrow \prod_{\mathfrak{p}\in \text{Minspec}(A)} {\rm GL}_{n}(A/\mathfrak{p}) , \]
		where $\text{Minspec}(A)$ is the set of all minimal prime ideals of $A$. As $A$ is noetherian, $\text{Minspec}(A)$ is a finite set. Thus, we assume further that $A$ is a pro-$p$ domain. Suppose that $G$ is a closed subgroup of ${\rm GL}_{n}(A)$ which contains a dense subset $D$ of torsions. By \cite[Proposition 3.2]{MR4976749}, there exists a positive integer $d$ and a family $\mathcal{N}$ of finite extensions of $\text{Frac}(A_{0})$ with degree less than or equal to $d$ such that the following homomorphism is injective:
		\[ {\rm GL}_{n}(A)\hookrightarrow \prod_{F\in \mathcal{N}}{\rm GL}_{n}(\mathcal{O}_{F}). \]
		Since $A$ is a pro-$p$ ring, ${\rm GL}_{n}(A)$ contains a closed pro-$p$ subgroup of finite index, and so does $G$. Let $P$ be a closed pro-$p$ subgroup $G$ of finite index. Then the image $\pi_{F}(P)$ of $P$ in ${\rm GL}_{n}(\mathcal{O}_{F})$ contains a dense subset of torsions for every $F\in \mathcal{N}$. Note that any pro-$p$ group with a solvable subgroup of finite index is actually solvable because any finite $p$-group is solvable. It follows that each $\pi_{F}(P)$ is solvable for every $F\in \mathcal{N}$. Moreover, $\pi_{F}(P)$ is solvable of derived length at most $ 2n$ for any $F\in \mathcal{N}$ by \cite[Corollary]{MR236272}. Since any product of solvable groups of a fixed derived length are solvable, $P$ is solvable. This completes the proof of our theorem.
	\end{proof}
	
	\subsection{Faltings' finiteness criteria}
	We will need the following important fact, which is due to Faltings.
	\begin{Lem}\label{Falting}
		Let $ K $ be a number field and $ S $ a finite set of primes of $ K $. Let $F$ be a non-Archimedean local field and $n$ be a positive integer. Then there exists a finite set $T$ of primes of $K$, depending on $n$ and disjoint from $S\cup S_{p}$, such that the isomorphism class of any semisimple continuous representation $\rho:G_{K,S\cup S_{p}}\to {\rm GL}_{n}(F)$ is determined uniquely by the finite collection $\{\operatorname{Tr}(\rho(\operatorname{Frob}_{v}))~|~v\in T\}$ where $\operatorname{Frob}_{v}$ denotes the Frobenius element at $v$ and $\operatorname{Tr}$ denotes the trace of square matrix.
	\end{Lem}
	\begin{proof}
		It is \cite[Section 2, Chapter V]{MR1175627}. See also \cite[Corollary 11]{MR4591759}.
	\end{proof}
	
	\subsection{Proof of Theorem \ref{local-globalprinciple}} 
	To prove Theorem \ref{local-globalprinciple}, we first establish the following result, which shows that Conjecture \ref{BUFM} depends only on the semisimplification of the representation.
	\begin{Prop}\label{semisimplesimplication}
		Let $ K $ be a number field and $ S $ a finite set of primes of $ K $ not containing any prime above $ p $. Let $F$ be a non-Archimedean local field with residue field of characteristic $p$ and let $\rho:G_{K,S}\to {\rm GL}_{n}(F)$ be a continuous representation. Then $\rho$ has finite image if and only if the semisimplification $\rho^{ss}$ of $\rho$ has finite image.
		
		Furthermore, there exists a finite set $T$ of primes of $K$, depending on $n$ and disjoint from $S\cup S_{p}$, with the following property: if $\rho(\operatorname{Frob}_{v})$ are unipotent for all $v\in T$, then $\rho$ has finite image where $\rho(\operatorname{Frob}_{v})$ denotes the Frobenius element corresponding to the prime $v$ in the group $\rho(G_{K})$.
	\end{Prop}
	\begin{proof}
		By \cite[Theorem 2.4]{MR4976749}, $\rho$ has finite image if and only if $\rho$ contains a solvable subgroup of finite index. By Chebotarev density theorem, the set of Frobenius elements $\operatorname{Frob}_{v}$ for all $v\notin S$ is dense in $G_{K,S}$. Then the same argument in the proof of Theorem \ref{niceobsection} shows that $\rho$ has finite image if and only if the eigenvalues of all $\rho(\operatorname{Frob}_{v}),v\notin S$, are roots of unity. It follows from Brauer–Nesbitt theorem (cf. \cite[Section 30.16]{MR0144979}) that the characteristic polynomials of $\rho(g)$ and $\rho^{ss}(g)$ are the same for all $g\in G_{K,S}$. Therefore, $\rho$ has finite image if and only if $\rho^{ss}$ has finite image. Moreover, our second claim follows from Lemma \ref{Falting}. This completes the proof of our proposition.
	\end{proof}
	
	We can now prove Theorem \ref{local-globalprinciple} as follows.
	\begin{proof}[Proof of Theorem \ref{local-globalprinciple}]
		Since $\rho(I_{v})$ is finite for any prime $v$ of $K$ above $p$, we may assume that $\rho$ factors through $G_{K,S}$ where $S$ is a finite set of primes of $K$ not containing any prime above after taking a finite extension of $K$. By our assumption and Chebotarev density theorem, the image $\text{Im}(\rho)$ of $\rho$ contains a dense subset of torsions. Then our first assertion follows from Theorem \ref{niceobsection}.
		
		For the second assertion, we can assume that $A=\mathcal{O}_{F}$ for some non-Archimedean local field $F$ with residue field of characteristic $p$ by \cite[Theorem 2.4 and 3.9]{MR4976749}. Then our assertion follows from Proposition \ref{semisimplesimplication}. The proof is finished.
	\end{proof}

	\section{A prime-to-adjoint principle}\label{section5}
	
	\subsection{Preliminaries} We begin by introducing some notations and conventions.
	\begin{Convention}
		\begin{enumerate}
			\item Let $K$ be a number field and $p$ be a prime number. If $S$ is a finite set of primes of $K$, then we denote by
			\[  V_{S}(K):=\{x\in K^{\times}~|~(x)=I^{p}~\text{as a fractional ideal of $K$},~x\in K_{v}^{p}~\text{for}~ v\in S\},\]
			where $K_{v}$ denotes the completion of $K$ at $v$. Moreover, we put
			\[ B_{S}(K):=V_{S}(K)/K^{\times p}. \]
			\item Let $G$ be a pro-$p$ group. We denote by $\overline{G}$ the Frattini quotient $G/\Phi(G)$ of $G$ where $\Phi(G)$ is the Frattini subgroup of $G$, i.e., the topological closure of $[G,G]G^{p}$ in $G$.
			\item If $K$ is a number field, then we denote by $ \operatorname{Cl}(K)$ the ideal class group of $K$ and $\operatorname{Cl}(K)[p]:=\{a\in \operatorname{Cl}(K)~|~pa=0\}$.
			\item If $H$ is a group and $\mathbb{F}$ is a field, then we denote by $\mathbb{F}[H]$ the group ring of $H$ over $\mathbb{F}$.
			\item If $A$ is a pro-$p$ ring with maximal ideal $\mathfrak{m}$, then we define $ {\rm GL}_{n}^{1}(A):=\ker({\rm GL}_{n}(A)\to {\rm GL}_{n}(A/\mathfrak{m}))$.
			\item Suppose that $\mathbb{F}$ is a finite field of characteristic $p$. If $H$ is a finite group of order prime to $p$ and $V_{1},V_{2}$ are two $\mathbb{F}[H]$-modules, then we say that $V_{1}$ is prime to $V_{2}$ if they have no irreducible subrepresentation in common.
		\end{enumerate}
	\end{Convention}
 
	\begin{Def}
	Let $\mathbb{F}$ be a finite field of characteristic $p$. Let $H$ be a subgroup of ${\rm GL}_{n}(\mathbb{F})$ whose order is prime to $p$ and $M$ be the corresponding adjoint $\mathbb{F}[H]$-module. A $\mathbb{F}[H]$-moduel $V$ is called \textit{prime-to-adjoint} if $V$ is prime to $M$ as $\mathbb{F}[H]$-modules.
	\end{Def}
	
	We will need the following lemma from group theory.
	\begin{Lem}\label{Bostonlemma}
		\begin{enumerate}
			\item Let $P$ be a pro-$p$ group. If $x_{1},\cdots,x_{d}\in P$ map to generators of $\overline{P}$, then they generate $P$. 
			\item Let $H$ be a finite group of order prime to $p$ and let $P$ be a pro-$p$ group. If $H$ acts on $P$ and $V$ is an $\mathbb{F}_{p}[H]$-submodule of $\overline{P}$, then one can find an $H$-invariant subgroup $B$ in $P$ with $\dim_{\mathbb{F}_{p}}V$ generators mapping onto $V$.
			\item Let $A$ be a pro-$p$ ring with residue field $\mathbb{F}$ and $X$ be a topologically finitely generated subgroup of ${\rm GL}_{n}^{1}(A)$. Let $H$ be a subgroup of ${\rm GL}_{n}(A)$ of order prime to $p$ that normalises $X$. Suppose that the $\mathbb{F}_{p}[H]$-module $\overline{X}$ is prime-to-adjoint. Then $X$ is trivial.
		\end{enumerate}
	\end{Lem}
	\begin{proof}
	These are \cite[Proposition 2.2, Proposition 2.3 and Lemma 2.5]{MR1079842}.
	\end{proof}

	Let $A$ be a pro-$p$ ring with maximal ideal $\mathfrak{m}$ and residue field $\mathbb F$. Let $H$ be a finite subgroup of ${\rm GL}_{n}(A)$ whose order is prime to $p$. Then we have an isomorphism $\pi:H\simeq \pi(H)$ of groups where $\pi$ is the group homomorphism from ${\rm GL}_{n}(A)$ to $ {\rm GL}_{n}(\mathbb{F})$ induced by the ring surjection $A\to A/\mathfrak{m}=\mathbb{F}$ because ${\rm GL}_{n}^{1}(A)$ is a pro-$p$ group. 
		
	\begin{Lem}\label{multiplicity-freecondition}
	With the notation as above. Suppose that the action of $\pi(H)=H$ on $V=\mathbb{F}^n$ is multiplicity-free, i.e., each irreducible constituent of the associated $\mathbb{F}[H]$-module appears with multiplicity exactly one. Then the fixed point subgroup
		\[
		({\rm GL}_n^1(A))^H
		\]
	   of ${\rm GL}_{n}^{1}(A)$ is abelian where $H$ acts on ${\rm GL}_{n}^{1}(A)$ by conjugation.
	\end{Lem}
	\begin{proof}
		Since the order of $H$ is prime to $p$, the $\mathbb F[H]$-module $V$ is semisimple. By the multiplicity-free assumption and Schur's lemma, the centralizer of $H$ in the algebra $M_n(\mathbb{F})$ of $n\times n$ matrices over $\mathbb{F}$ is commutative.
		
		Since $	{\rm GL}_{n}^1(A) =1+M_n(\mathfrak m)$, any element of $({\rm GL}_n^1(A))^H$ can be written as $1+X$ with
		$X\in M_n(\mathfrak m)$ fixed under conjugation by $H$. Since $A$ is Noetherian and $\mathfrak m$-adically complete, we have $\bigcap_{k\ge 1}\mathfrak m^k=0$. Thus, it suffices to show that
		\[ [X,Y]:=XY-XY=0\pmod{\mathfrak{m}^{k}} \]
		  for all $X,Y\in M_{n}(\mathfrak{m})$ fixed under the conjugation by $H$ and all $k$. It is trivial that the claim holds for $k=1$. We proceed by induction on $k$. Assume the claim holds for some $k\ge 1$. Then we have 
		\[ [X,Y]=tM \pmod{\mathfrak{m}^{k+1}} \]
		where $t\in \mathfrak{m}^{k}$ and $M\in M_{n}(A)$. Since the centralizer of $H$ in $M_n(\mathbb{F})$ is commutative, we have $ M=0\pmod{\mathfrak{m}}$. It follows that $tM=0 \pmod{\mathfrak{m}^{k+1}}$. This finishes the proof of our claim.
	\end{proof}
	
	We need the following simple algebraic fact.
		\begin{Lem}\label{simplelemma}
		Let $H$ be a finite group whose order is prime to $p$. Suppose that the following sequences of $\mathbb{F}_{p}[H]$-modules
		\[ 0\to M_{1}\to M_{2}\to M_{3}\to M_{4}\to M_{5}\to 0, \]
		\[ 0\to N_{1}\to M_{2}\to N_{2}\to 0,  \]
		are exact where $N_{2}\cong M_{5}$. Then we have an isomorphism of $\mathbb{F}_{p}[H]$-modules
		\[ M_{4}\cong M_{1}/(M_{1}\cap N_{1})\oplus \operatorname{Coker}(N_{1}\to M_{3}). \]
	\end{Lem}
	\begin{proof}
		Since the order of $H$ is prime to $p$, any short exact sequence of $\mathbb{F}_{p}[H]$-modules splits. The claim follows.
	\end{proof}
	
The following theorem on the structure of the Galois group is needed.
   	\begin{Prop}\label{structureofGLs}
   	Let $L$ be a finite Galois extension of $K$ with Galois group $H$ whose order is to $p$. Let $S$ be a finite set of primes of $K$ not containing any prime above $p$ and $S(L)$ be the set of primes of $L$ above $S$. Let $\overline{E}(L)$ and $\overline{E}_{v}(L)$ be the global and local units modulo $p$-powers of $L$ and $L_{v}$ respectively where $v\in S$. Then we have an isomorphism of $\mathbb{F}_{p}[H]$-modules
   	\[ \overline{G_{L,S(L)}(p)}\cong B_{S(L)}(L)/(B_{S(L)}(L)\cap \overline{E}(L))\oplus  \operatorname{Coker}(\overline{E}(L)\to \bigoplus_{v\in S(L)}\mu_{p}(L_{v})),\]
   	where $\mu_{p}(L_{v})$ denotes the group of $p$-th roots of unity in $L_{v}$.
   \end{Prop}
   \begin{proof}
   	First note that we have an exact sequence of $\mathbb{F}_{p}[H]$-modules
   	\[ 	0\to \overline{E}(L)\to B_{\emptyset}(L)\to \operatorname{Cl}(L)[p]\to 0.  \]
   	By \cite[Section 11.3]{MR1930372}, we have the following exact sequence of $\mathbb{F}_{p}[H]$-modules
   	\[ 0\to B_{S(L)}(L)\to B_{\emptyset}(L)\to \bigoplus_{v\in S(L)}\overline{E}_{v}(L)\to \overline{G_{L,S}(p)}\to \operatorname{Cl}(L)/p\to 0, \]
   	By \cite[Lemma 3.3]{MR1188818}, we have an isomorphism of $\mathbb{F}_{p}[H]$-modules
   	\[ \operatorname{Cl}(L)[p]\cong \operatorname{Cl}(L)/p .\]
   	Since the prime $v\in S(L)$ does not lie above $p$, we also have an isomorphism of $\mathbb{F}_{p}[H]$-modules
   	\[ \bigoplus_{v\in S(L)}\overline{E}_{v}(L)\cong \bigoplus_{v\in S(L)}\mu_{p}(L_{v}). \]
   	Then our claim follows from Lemma \ref{simplelemma}.
   \end{proof}	
   
     \subsection{A prime-to-adjoint principle}
   \begin{Convention}
   	Let $\mathbb{F}$ be a finite field of characteristic $p$. If $H$ is a finite group of order prime to $p$ and $V$ is a $ \mathbb{F}[H]$-module, then we denote by $V^{\dagger}$ the quotient $\mathbb{F}[H]$-module $V/V^{H}$.
   \end{Convention}
    In this section, our main result is the following.
	\begin{Thm}\label{ptatheorem}
		Let $K$ be a number field and $S$ a finite set of primes of $K$ not containing any prime above $p$. Let $\rho:G_{K,S}\to {\rm GL}_{n}(A)$ be a continuous homomorphism where $ A $ is a pro-$p$ ring with residue field $\mathbb{F}$. Assume the following conditions hold.
		\begin{enumerate}
			\item The image of $\text{Im}(\overline{\rho})$ of $\overline{\rho}$ contains a group $H$ whose order prime to $p$.
			\item The representation $\overline{\rho}|_{H}$ is multiplicity-free.
			\item The adjoint representation $\operatorname{ad}(\overline{\rho}|_{H})$ of $\rho|_{H}$ is prime to
			\begin{equation}\label{primeto}
			 \left(  B_{S(L)}(L)/(B_{S(L)}(L)\cap \overline{E}(L))\oplus  \operatorname{Coker}(\overline{E}(L)\to \bigoplus_{v\in S(L)}\mu_{p}(L_{v}))\right)^{\dagger}
			\end{equation}
		  where $L$ denotes the splitting field of $\overline{\rho}$, i.e., the Galois extension $L$ of $K$ such that $\operatorname{Gal}(L/K)=\text{Im}(\overline{\rho})$.
		\end{enumerate}
		Then $\rho$ has finite image.
	\end{Thm}
	\begin{proof}
		Let $L$ be the splitting field of $\overline{\rho}$. After taking a finite extension of $K$, we may assume that $\text{Im}(\overline{\rho})=H$. Let $S(L)$ be the set of primes of $K$ above $S$. Then it suffices to show that
		\[\rho|_{\ker(\overline{\rho})}\in \operatorname{Hom}_{H}(G_{L,S(L)}(p),{\rm GL}_{n}^{1}(A))  \]
		 has finite image. By \cite[Theorem 11.5]{MR1930372}, the pro-$p$ group $ G_{L,S(L)}(p)$ is finitely generated. Let $x_{1},\cdots,x_{d}\in G_{L,S(L)}(p) $ map to generators $\bar{x}_{1},\cdots,\bar{x}_{d}$ of $\overline{G_{L,S(L)}(p)}$. By our assumption (3) and Proposition \ref{structureofGLs}, the set of generators $\{\bar{x}_{1},\cdots,\bar{x}_{d}\}$ consists of two disjoint subsets: the fixed points and those that are prime-to-adjoint.
		By Lemma \ref{Bostonlemma}(3), genrators that are prime-to-adjoint go to trivial. It follows from Lemma \ref{Bostonlemma} that the image of $\rho|_{\ker(\overline{\rho})}$ is contained in $({\rm GL}_n^1(A))^H$. By our assumption (2) and Lemma \ref{multiplicity-freecondition}, the image of $\rho|_{\ker(\overline{\rho})}$ is abelian. It follows that the image of $\rho$ contains an open abelian subgroup. Finally, our claim follows from \cite[Theorem 2.4]{MR4976749}.
	\end{proof}
	
    Finally, we prove Theorem \ref{primetoadjointprincple} as follows.
    
		\begin{proof}[Proof of Theorem \ref{primetoadjointprincple}]
	If $S=\emptyset$, then we have $\overline{G_{L,\emptyset}(p)}=\operatorname{Cl}(L)/p$ and (\ref{primeto}) in Theorem \ref{ptatheorem} is just $\ker(N:\operatorname{Cl}(L)/p\to \operatorname{Cl}(K')/p)$. Then our claim follows from Theorem \ref{ptatheorem}.
	\end{proof}

	\section*{Acknowledgements}
     The author is grateful to Professor Mikhail Ershov for helpful correspondence regarding Golod-Shafarevich groups.

	\bibliographystyle{plain}
	\bibliography{FM}

@Article{MR1981910,
  author     = {Kisin, M. and Wortmann, S.},
  journal    = {Math. Res. Lett.},
  title      = {A note on {A}rtin motives},
  year       = {2003},
  issn       = {1073-2780},
  number     = {2-3},
  pages      = {375--389},
  volume     = {10},
  doi        = {10.4310/MRL.2003.v10.n3.a7},
  fjournal   = {Mathematical Research Letters},
  mrclass    = {14F30 (14F42 14G25)},
  mrnumber   = {1981910},
  mrreviewer = {Tam\'{a}s Szamuely},
  url        = {https://doi.org/10.4310/MRL.2003.v10.n3.a7},
}

@InCollection{MR1363495,
  author     = {Fontaine, J.-M. and Mazur, B.},
  booktitle  = {Elliptic curves, modular forms, \& {F}ermat's last theorem ({H}ong {K}ong, 1993)},
  publisher  = {Int. Press, Cambridge, MA},
  title      = {Geometric {G}alois representations},
  year       = {1995},
  pages      = {41--78},
  series     = {Ser. Number Theory, I},
  mrclass    = {11F80 (11G35 11R32)},
  mrnumber   = {1363495},
  mrreviewer = {Ian Kiming},
}

@Book{MR2392026,
  author    = {Neukirch, J. and Schmidt, A. and Wingberg, K.},
  publisher = {Springer-Verlag, Berlin},
  title     = {Cohomology of number fields},
  year      = {2008},
  edition   = {Second},
  isbn      = {978-3-540-37888-4},
  series    = {Grundlehren der mathematischen Wissenschaften [Fundamental Principles of Mathematical Sciences]},
  volume    = {323},
  doi       = {10.1007/978-3-540-37889-1},
  mrclass   = {11R34 (11-02 11G45 11R23 11S20 11S25 11S31 12G05)},
  mrnumber  = {2392026},
  pages     = {xvi+825},
  url       = {https://doi.org/10.1007/978-3-540-37889-1},
}

@Book{MR554237,
  author    = {Serre, J.-P.},
  publisher = {Springer-Verlag, New York-Berlin},
  title     = {Local fields},
  year      = {1979},
  isbn      = {0-387-90424-7},
  note      = {Translated from the French by Marvin Jay Greenberg},
  series    = {Graduate Texts in Mathematics},
  volume    = {67},
  mrclass   = {12Bxx},
  mrnumber  = {554237},
  pages     = {viii+241},
}

@Book{MR1930372,
  author    = {Koch, H.},
  publisher = {Springer-Verlag, Berlin},
  title     = {Galois theory of {$p$}-extensions},
  year      = {2002},
  isbn      = {3-540-43629-4},
  note      = {With a foreword by I. R. Shafarevich, Translated from the 1970 German original by Franz Lemmermeyer, With a postscript by the author and Lemmermeyer},
  series    = {Springer Monographs in Mathematics},
  doi       = {10.1007/978-3-662-04967-9},
  mrclass   = {11S25 (11R32 11R34 11S20)},
  mrnumber  = {1930372},
  pages     = {xiv+190},
  url       = {https://doi.org/10.1007/978-3-662-04967-9},
}

@Article{MR1681626,
  author     = {Boston, N.},
  journal    = {J. Number Theory},
  title      = {Some cases of the {F}ontaine-{M}azur conjecture. {II}},
  year       = {1999},
  issn       = {0022-314X},
  number     = {2},
  pages      = {161--169},
  volume     = {75},
  doi        = {10.1006/jnth.1998.2337},
  fjournal   = {Journal of Number Theory},
  mrclass    = {11R32 (11R37)},
  mrnumber   = {1681626},
  mrreviewer = {Franz Lemmermeyer},
  url        = {https://doi.org/10.1006/jnth.1998.2337},
}

@Book{MR1720368,
  author     = {Dixon, J. D. and du Sautoy, M. P. F. and Mann, A. and Segal, D.},
  publisher  = {Cambridge University Press, Cambridge},
  title      = {Analytic pro-{$p$} groups},
  year       = {1999},
  edition    = {Second},
  isbn       = {0-521-65011-9},
  series     = {Cambridge Studies in Advanced Mathematics},
  volume     = {61},
  doi        = {10.1017/CBO9780511470882},
  mrclass    = {20E18 (20G30)},
  mrnumber   = {1720368},
  mrreviewer = {Alexander Lubotzky},
  pages      = {xviii+368},
  url        = {https://doi.org/10.1017/CBO9780511470882},
}

@Article{MR1432356,
  author     = {Hajir, F.},
  journal    = {J. Algebra},
  title      = {On the growth of {$p$}-class groups in {$p$}-class field towers},
  year       = {1997},
  issn       = {0021-8693},
  number     = {1},
  pages      = {256--271},
  volume     = {188},
  doi        = {10.1006/jabr.1996.6849},
  fjournal   = {Journal of Algebra},
  mrclass    = {11R37 (20E18)},
  mrnumber   = {1432356},
  mrreviewer = {Marcus du Sautoy},
  url        = {https://doi.org/10.1006/jabr.1996.6849},
}

@Article{MR1079842,
  author     = {Boston, N.},
  journal    = {Invent. Math.},
  title      = {Explicit deformation of {G}alois representations},
  year       = {1991},
  issn       = {0020-9910},
  number     = {1},
  pages      = {181--196},
  volume     = {103},
  doi        = {10.1007/BF01239511},
  fjournal   = {Inventiones Mathematicae},
  mrclass    = {11F80 (11G05)},
  mrnumber   = {1079842},
  mrreviewer = {Bas Edixhoven},
  url        = {https://doi.org/10.1007/BF01239511},
}

@Article{MR0161852,
  author     = {Golod, E. S. and \v{S}afarevi\v{c}, I. R.},
  journal    = {Izv. Akad. Nauk SSSR Ser. Mat.},
  title      = {On the class field tower},
  year       = {1964},
  issn       = {0373-2436},
  pages      = {261--272},
  volume     = {28},
  fjournal   = {Izvestiya Akademii Nauk SSSR. Seriya Matematicheskaya},
  mrclass    = {10.68},
  mrnumber   = {0161852},
  mrreviewer = {E. Inaba},
}

@Article{MR4308183,
  author     = {Hajir, F. and Maire, C. and Ramakrishna, R.},
  journal    = {Ann. Math. Qu\'{e}.},
  title      = {Cutting towers of number fields},
  year       = {2021},
  issn       = {2195-4755},
  number     = {2},
  pages      = {321--345},
  volume     = {45},
  doi        = {10.1007/s40316-021-00156-8},
  fjournal   = {Annales Math\'{e}matiques du Qu\'{e}bec},
  mrclass    = {11R29 (11R21 11R37)},
  mrnumber   = {4308183},
  mrreviewer = {Abdelmalek Azizi},
  url        = {https://doi.org/10.1007/s40316-021-00156-8},
}

@Article{MR1935507,
  author     = {Jaikin-Zapirain, A.},
  journal    = {J. Algebra},
  title      = {On linear just infinite pro-{$p$} groups},
  year       = {2002},
  issn       = {0021-8693},
  number     = {2},
  pages      = {392--404},
  volume     = {255},
  doi        = {10.1016/S0021-8693(02)00024-8},
  fjournal   = {Journal of Algebra},
  mrclass    = {20E18 (20G25)},
  mrnumber   = {1935507},
  mrreviewer = {Martyn Quick},
  url        = {https://doi.org/10.1016/S0021-8693(02)00024-8},
}

@InCollection{MR1765119,
  author    = {Grigorchuk, R. I.},
  booktitle = {New horizons in pro-{$p$} groups},
  publisher = {Birkh\"{a}user Boston, Boston, MA},
  title     = {Just infinite branch groups},
  year      = {2000},
  pages     = {121--179},
  series    = {Progr. Math.},
  volume    = {184},
  mrclass   = {20F14 (20E08 20E18)},
  mrnumber  = {1765119},
}

@Book{MR0144979,
  author     = {Curtis, C. W. and Reiner, I.},
  publisher  = {Interscience Publishers (a division of John Wiley \& Sons, Inc.), New York-London},
  title      = {Representation theory of finite groups and associative algebras},
  year       = {1962},
  series     = {Pure and Applied Mathematics, Vol. XI},
  mrclass    = {20.80},
  mrnumber   = {0144979},
  mrreviewer = {W. E. Jenner},
  pages      = {xiv+685},
}

@Article{MR1264349,
  author     = {Lubotzky, A. and Shalev, A.},
  journal    = {Israel J. Math.},
  title      = {On some {$\Lambda$}-analytic pro-{$p$} groups},
  year       = {1994},
  issn       = {0021-2172,1565-8511},
  number     = {1-3},
  pages      = {307--337},
  volume     = {85},
  doi        = {10.1007/BF02758646},
  fjournal   = {Israel Journal of Mathematics},
  mrclass    = {20E18 (22E20)},
  mrnumber   = {1264349},
  mrreviewer = {Marcus\ du Sautoy},
  url        = {https://doi.org/10.1007/BF02758646},
}

@Article{MR1684856,
  author     = {Barnea, Y. and Larsen, M.},
  journal    = {J. Algebra},
  title      = {A non-abelian free pro-{$p$} group is not linear over a local field},
  year       = {1999},
  issn       = {0021-8693,1090-266X},
  number     = {1},
  pages      = {338--341},
  volume     = {214},
  doi        = {10.1006/jabr.1998.7682},
  fjournal   = {Journal of Algebra},
  mrclass    = {20E06 (20G25)},
  mrnumber   = {1684856},
  mrreviewer = {Javier\ Otal},
  url        = {https://doi.org/10.1006/jabr.1998.7682},
}

@Article{MR1637068,
  author     = {Pink, R.},
  journal    = {J. Algebra},
  title      = {Compact subgroups of linear algebraic groups},
  year       = {1998},
  issn       = {0021-8693,1090-266X},
  number     = {2},
  pages      = {438--504},
  volume     = {206},
  doi        = {10.1006/jabr.1998.7439},
  fjournal   = {Journal of Algebra},
  mrclass    = {20G25 (22C05)},
  mrnumber   = {1637068},
  mrreviewer = {Herbert\ Abels},
  url        = {https://doi.org/10.1006/jabr.1998.7439},
}

@Article{MR2363421,
  author     = {Jaikin-Zapirain, A. and Klopsch, B.},
  journal    = {J. Lond. Math. Soc. (2)},
  title      = {Analytic groups over general pro-{$p$} domains},
  year       = {2007},
  issn       = {0024-6107,1469-7750},
  number     = {2},
  pages      = {365--383},
  volume     = {76},
  doi        = {10.1112/jlms/jdm055},
  fjournal   = {Journal of the London Mathematical Society. Second Series},
  mrclass    = {20E18 (20G25 22E20)},
  mrnumber   = {2363421},
  mrreviewer = {Nikolay\ V.\ Nikolov},
  url        = {https://doi.org/10.1112/jlms/jdm055},
}

@Article{MR4976749,
  author   = {Luo, Y.},
  journal  = {J. Number Theory},
  title    = {Remarks on the {B}oston's unramified {F}ontaine-{M}azur conjecture},
  year     = {2026},
  issn     = {0022-314X,1096-1658},
  pages    = {96--109},
  volume   = {281},
  doi      = {10.1016/j.jnt.2025.09.019},
  fjournal = {Journal of Number Theory},
  mrclass  = {11F80 (20E18)},
  mrnumber = {4976749},
  url      = {https://doi.org/10.1016/j.jnt.2025.09.019},
}

@InCollection{MR2148459,
  author     = {Boston, N.},
  booktitle  = {Progress in {G}alois theory},
  publisher  = {Springer, New York},
  title      = {Reducing the {F}ontaine-{M}azur conjecture to group theory},
  year       = {2005},
  isbn       = {0-387-23533-7},
  pages      = {39--50},
  series     = {Dev. Math.},
  volume     = {12},
  doi        = {10.1007/0-387-23534-5\_3},
  mrclass    = {14G32 (11R32 12F10 20E18)},
  mrnumber   = {2148459},
  mrreviewer = {Rachel\ J.\ Pries},
  url        = {https://doi.org/10.1007/0-387-23534-5_3},
}

@Article{MR236272,
  author     = {Dixon, J. D.},
  journal    = {Math. Z.},
  title      = {The solvable length of a solvable linear group},
  year       = {1968},
  issn       = {0025-5874,1432-1823},
  pages      = {151--158},
  volume     = {107},
  doi        = {10.1007/BF01111027},
  fjournal   = {Mathematische Zeitschrift},
  mrclass    = {20.75},
  mrnumber   = {236272},
  mrreviewer = {J.\ S.\ Rose},
  url        = {https://doi.org/10.1007/BF01111027},
}

@Article{MR4591759,
  author     = {Das, P. and Rajan, C. S.},
  journal    = {Proc. Amer. Math. Soc.},
  title      = {Finiteness theorems for potentially equivalent {G}alois representations: extension of {F}altings' finiteness criteria},
  year       = {2023},
  issn       = {0002-9939,1088-6826},
  number     = {8},
  pages      = {3189--3200},
  volume     = {151},
  doi        = {10.1090/proc/15856},
  fjournal   = {Proceedings of the American Mathematical Society},
  mrclass    = {11F80 (11G05 11G10)},
  mrnumber   = {4591759},
  mrreviewer = {Shiang\ Tang},
  url        = {https://doi.org/10.1090/proc/15856},
}

@Article{MR2332053,
  author     = {Boston, N.},
  journal    = {J. Th\'eor. Nombres Bordeaux},
  title      = {Galois groups of tamely ramified {$p$}-extensions},
  year       = {2007},
  issn       = {1246-7405,2118-8572},
  number     = {1},
  pages      = {59--70},
  volume     = {19},
  doi        = {10.5802/jtnb.573},
  fjournal   = {Journal de Th\'eorie des Nombres de Bordeaux},
  mrclass    = {11R32},
  mrnumber   = {2332053},
  mrreviewer = {Hiroaki\ Nakamura},
  url        = {https://doi.org/10.5802/jtnb.573},
}

@Article{MR707163,
  author     = {Lubotzky, A.},
  journal    = {Ann. of Math. (2)},
  title      = {Group presentation, {$p$}-adic analytic groups and lattices in {${\rm SL}\sb{2}({\bf C})$}},
  year       = {1983},
  issn       = {0003-486X,1939-8980},
  number     = {1},
  pages      = {115--130},
  volume     = {118},
  doi        = {10.2307/2006956},
  fjournal   = {Annals of Mathematics. Second Series},
  mrclass    = {22E40 (20F05)},
  mrnumber   = {707163},
  mrreviewer = {S.\ P.\ Demushkin},
  url        = {https://doi.org/10.2307/2006956},
}

@Article{MR172892,
  author     = {Vinberg, \`E.\ B.},
  journal    = {Izv. Akad. Nauk SSSR Ser. Mat.},
  title      = {On the theorem concerning the infinite-dimensionality of an associative algebra},
  year       = {1965},
  issn       = {0373-2436},
  pages      = {209--214},
  volume     = {29},
  fjournal   = {Izvestiya Akademii Nauk SSSR. Seriya Matematicheskaya},
  mrclass    = {10.68 (16.10)},
  mrnumber   = {172892},
  mrreviewer = {E.\ Inaba},
}

@Book{MR316588,
  author     = {Satake, I.},
  publisher  = {Marcel Dekker, Inc., New York},
  title      = {Classification theory of semi-simple algebraic groups},
  year       = {1971},
  note       = {With an appendix by M. Sugiura, Notes prepared by Doris Schattschneider},
  series     = {Lecture Notes in Pure and Applied Mathematics},
  volume     = {3},
  mrclass    = {20G15},
  mrnumber   = {316588},
  mrreviewer = {T.\ Ono},
  pages      = {viii+149},
}

@Book{MR1175627,
  author     = {Faltings, G. and W\"ustholz, G. and Grunewald, F. and Schappacher, N. and Stuhler, U.},
  publisher  = {Friedr. Vieweg \& Sohn, Braunschweig},
  title      = {Rational points},
  year       = {1992},
  edition    = {Third},
  isbn       = {3-528-28593-1},
  note       = {Papers from the seminar held at the Max-Planck-Institut f\"ur Mathematik, Bonn/Wuppertal, 1983/1984, With an appendix by W\"ustholz},
  series     = {Aspects of Mathematics},
  volume     = {E6},
  doi        = {10.1007/978-3-322-80340-5},
  mrclass    = {11G35 (11D41 14C17 14G05 14G40)},
  mrnumber   = {1175627},
  mrreviewer = {Kenneth\ A.\ Ribet},
  pages      = {x+311},
  url        = {https://doi.org/10.1007/978-3-322-80340-5},
}

@InCollection{MR1754662,
  author     = {Grigorchuk, R. I. and Herfort, W. N. and Zalesskii, P. A.},
  booktitle  = {Algebra ({M}oscow, 1998)},
  publisher  = {de Gruyter, Berlin},
  title      = {The profinite completion of certain torsion {$p$}-groups},
  year       = {2000},
  isbn       = {3-11-016399-3},
  pages      = {113--123},
  mrclass    = {20E18 (20E08)},
  mrnumber   = {1754662},
  mrreviewer = {Thomas\ Weigel},
}

@Article{Luo22012026,
  author    = {Y. Luo},
  journal   = {Communications in Algebra},
  title     = {The topological Tits alternative for linear groups over commutative profinite rings},
  year      = {2026},
  number    = {online},
  doi       = {10.1080/00927872.2025.2608664},
  eprint    = {https://doi.org/10.1080/00927872.2025.2608664},
  publisher = {Taylor \& Francis},
  url       = {https://doi.org/10.1080/00927872.2025.2608664},
}

@Article{MR2949205,
  author     = {Ershov, M.},
  journal    = {Internat. J. Algebra Comput.},
  title      = {Golod-{S}hafarevich groups: a survey},
  year       = {2012},
  issn       = {0218-1967,1793-6500},
  number     = {5},
  pages      = {1230001, 68},
  volume     = {22},
  doi        = {10.1142/S0218196712300010},
  fjournal   = {International Journal of Algebra and Computation},
  mrclass    = {20F05 (16W50 20E18 20F50 20F69)},
  mrnumber   = {2949205},
  mrreviewer = {Andrei\ Jaikin-Zapirain},
  url        = {https://doi.org/10.1142/S0218196712300010},
}

@InCollection{MR1765122,
  author     = {Zelmanov, E.},
  booktitle  = {New horizons in pro-{$p$} groups},
  publisher  = {Birkh\"auser Boston, Boston, MA},
  title      = {On groups satisfying the {G}olod-{S}hafarevich condition},
  year       = {2000},
  isbn       = {0-8176-4171-8},
  pages      = {223--232},
  series     = {Progr. Math.},
  volume     = {184},
  mrclass    = {20E18 (17B01)},
  mrnumber   = {1765122},
  mrreviewer = {Alexander\ Lubotzky},
}

@Article{MR1109625,
  author     = {Wilson, J. S.},
  journal    = {Invent. Math.},
  title      = {Finite presentations of pro-{$p$} groups and discrete groups},
  year       = {1991},
  issn       = {0020-9910,1432-1297},
  number     = {1},
  pages      = {177--183},
  volume     = {105},
  doi        = {10.1007/BF01232262},
  fjournal   = {Inventiones Mathematicae},
  mrclass    = {20F05 (20E18)},
  mrnumber   = {1109625},
  mrreviewer = {Marcus\ du Sautoy},
  url        = {https://doi.org/10.1007/BF01232262},
}

@InCollection{MR224710,
  author     = {Tits, J.},
  booktitle  = {Algebraic {G}roups and {D}iscontinuous {S}ubgroups ({P}roc. {S}ympos. {P}ure {M}ath., {B}oulder, {C}olo., 1965)},
  publisher  = {Amer. Math. Soc., Providence, RI},
  title      = {Classification of algebraic semisimple groups},
  year       = {1966},
  pages      = {33--62},
  mrclass    = {20.27},
  mrnumber   = {224710},
  mrreviewer = {R.\ Steinberg},
}

@Article{MR1188818,
  author     = {Boston, N. and Ullom, S. V.},
  journal    = {Math. Proc. Cambridge Philos. Soc.},
  title      = {Representations related to {CM} elliptic curves},
  year       = {1993},
  issn       = {0305-0041,1469-8064},
  number     = {1},
  pages      = {71--85},
  volume     = {113},
  doi        = {10.1017/S0305004100075770},
  fjournal   = {Mathematical Proceedings of the Cambridge Philosophical Society},
  mrclass    = {11G05 (11F80)},
  mrnumber   = {1188818},
  mrreviewer = {Matthias\ Flach},
  url        = {https://doi.org/10.1017/S0305004100075770},
}

@Book{MR1978431,
  author     = {Lubotzky, A. and Segal, D.},
  publisher  = {Birkh\"auser Verlag, Basel},
  title      = {Subgroup growth},
  year       = {2003},
  isbn       = {3-7643-6989-2},
  series     = {Progress in Mathematics},
  volume     = {212},
  doi        = {10.1007/978-3-0348-8965-0},
  mrclass    = {20E07 (20E18 20E26 20F69)},
  mrnumber   = {1978431},
  mrreviewer = {Avinoam\ Mann},
  pages      = {xxii+453},
  url        = {https://doi.org/10.1007/978-3-0348-8965-0},
}

\end{document}